\documentclass{article}

\usepackage{amssymb}
\usepackage[english]{babel}
\usepackage[latin1]{inputenc}
\usepackage{graphicx}
\usepackage{epsfig}
\usepackage{amsmath}
\usepackage{amsthm}

 \newtheorem{Theo}{Theorem}[section]
 \newtheorem{Co}[Theo]{Corollary}
 \newtheorem{Le}[Theo]{Lemma}
 \newtheorem{Pro}[Theo]{Proposition}
\theoremstyle{definition}
 
 \theoremstyle{remark}
 
 \numberwithin{equation}{section}
\newtheorem*{Con}{Conjecture}

\newcommand{\avoir}{generic }
\newcommand{\alignment}{alignment }
\newcommand{\grad}{\nabla}
\newcommand{\imprimitive}{$2$-imprimitive }
\newcommand{\imprimitively}{$2$-imprimitively }
\newcommand{\imprimitivity}{$2$-imprimitivity }
\newcommand{\mm}{\mathcal{M}}

\newcommand{\pp}{\mathcal{P}}

\newcommand{\lk}{\mathcal{L}}
\newcommand{\hh}{\mathcal{H}}

\newcommand{\cc}{\mathcal{C}}
\newcommand{\ii}{\mathcal{I}}
\newcommand{\vv}{\mathcal{V}}
\newcommand{\uu}{\mathcal{U}}

\newcommand{\g}{\mathfrak{g}}
\newcommand{\avoirity}{genericity }

\begin{document}

\title{Turbiner's Conjecture in Three Dimensions}

\author{Mélisande Fortin Boisvert \thanks{Department of Mathematics and Statistics, McGill University, Montréal, Québec, H3A~2K6,
 Canada. E-mail: boisvert@math.mcgill.ca}}







\date{}




\maketitle

\begin{abstract}
We prove a modified version of Turbiner's conjecture in three
dimensions and we give a counter-example to the original conjecture.
The Lie algebraic Schrödinger operators corresponding to flat
metrics of a certain restricted type are shown to separate partially
in either Cartesian, cylindrical or spherical coordinates.
\end{abstract}

           \section{Introduction}
            The aim of this article is to extend results related to separation of
            variables
          for flat Lie algebraic Schrödinger operators. Originally, in \cite{Turbiner}, Alexander
          Turbiner conjectured the following:

\begin{Con}(Turbiner) In $\mathbb{R}^2$ there exist no quasi-exactly-solvable or exactly-solvable
problems containing the Laplace-Beltrami operator with flat-space
metric tensor, which are characterized by non-separable
variables.\end{Con} This conjecture was reformulated in more
geometrical terms by Rob Milson. The conjecture, which the work in
this paper is based on, now reads as follows.

           \begin{Con}(Turbiner, second version) Let H be a Lie algebraic Schrödinger operator defined
on a 2-dimensional manifold. If the symbol of H engenders a
Euclidean geometry, i.e. if the corresponding Gaussian curvature is
zero, then the spectral equation $H\psi = E\psi$ can be solved by a
separation of variables.\end{Con}

   This conjecture is false in general. A counter-example is given by Rob Milson
  in \cite{MilsonThesis} and
  \cite{Imprimitively}, together with a proof of a modified version of the conjecture. By adding
  two extra assumptions, namely an imprimitive action and a compactness requirement, one
  can prove that the spectral equation can be solved by separation of variables.
Furthermore, the imprimitivity hypothesis implies even more than
expected: separation will occur in either a Cartesian
or polar coordinate system.\\

  In this paper, it is shown that, in three dimensions, the original conjecture is also false and
  a proof of a modified version of the $3$D Turbiner's conjecture is given. Again the two preceding requirements
   are indispensible, and
  a third condition, related to the contravariant metric, will be
imposed. Once again, the imprimitivity of the action will
  ensure that separation, here only partial, will occur in either a Cartesian, cylindrical or
  spherical coordinate system.\\

 The proofs of both modified versions of the conjecture are based on the the following ideas.
 First, the imprimitive action induces an invariant foliation
                                  $\Lambda$, for which the leaves
                                  will be denoted by
                                  $\{\lambda=\textrm{constant}\}$. The Schrödinger operator
                                  $\hh$ is Lie algebraic, thus, it
                                  is an element of the enveloping algebra of a Lie
                                  algebra of first order differential
                                  operators.  When applied to
                                  $\lambda$, the elements of the generating Lie algebra
                                  must give back functions of
                                  $\lambda$. The
                                  operator $\hh$  will enjoy the same property, that is $\hh(\lambda)=f(\lambda)$.
                                  Combining the fact that the operator is Lie algebraic with the
                                  imprimitivity of the action, one can prove that the leaves of the
                                  foliation which is
                                  perpendicular to $\Lambda$ are
                                  necessarily geodesics.
                                  Then, one of the key tool, the Tiling Theorem, gives a
                                  global map from a Euclidean
                                  space to our manifold. Thus
                                  by pulling back the operator, the
                                  leaves of the perpendicular foliation are straight
                                  lines. Note that this theorem follows from a intermediate one: the Trapping
                                  Theorem.\\

                                 In this setting, one can
                                  show that the invariant leaves can only be
                                  prescribed curves or surfaces. In
                                  two dimensions the curves need to be
                                  straight lines or concentric circles, while in three
                                  dimensions the surfaces have to be
                                  planes, cylinders or spheres.
                                  In this context, $\lambda$ will be
                                  either a Cartesian or a radial
                                  coordinate. Finally, using the appropriate coordinate system, one checks that the equation $\hh \psi=E \psi$
                                  separates with respect to
                                  the coordinate $\lambda$.\\

 Despite the fact that the path followed to prove the
  modified $3$D Turbiner's conjecture is similar to the one given in
  \cite{Imprimitively}, there are several important issues, which were absent in the two
  dimensional case, and which appear in our study. We had first to select the appropriate
  flatness criteria, while in two dimensions there is no
    such choice. In order to prove the $3$D-Trapping Theorem, one has to
    impose the diagonal
    terms of the Ricci curvature tensor to be zero. For the $3$D-Tiling
    Theorem, it is the Riemannian curvature tensor that needs to vanish.
   For this proof of the $3$D-Trapping and Tiling Theorems, we have to assume
   that either the metric can be diagonalised, either $\mathbf{M}$ is a transverse, type changing
   manifold. For the first case, to conclude both theorems, an extra requirement: the
\avoirity of the contravariant metric needs to be added. This
requirement is related to the non-invertible factors of its
components and will be defined later. For the transverse type
changing manifold, we will see that such
metric can always be diagonalised and is necessarily \avoir. \\

         The determination of the possible foliations required a different approach. While in
         two dimensions
    one could only consider the possible
foliations by straight lines to conclude, in three dimensions, one
has to keep in mind the entire picture of the two perpendicular
foliations in order to determine the three types of leaves. The
arguments are not sophisticated but the proof is long enough for us
to devote an entire section to it.
 Another distinction with the two dimensional case is
the fact that separation of variables is only partial. Indeed, as in
the work of Rob Milson, one can isolate one variable, $\lambda$, but
we are left,
in three dimensions, with two variables for which nothing can be said.\\

 In section~$2$, we briefly describe the context of Turbiner's conjecture. We define all the notions employed
   in this paper and we give an example that illustrate how the imprimitivity of the action
induces separation of variables.
   Section~$3$ fills in the gaps needed to generalize the proof of both $3$D-Tiling and
   $3$D-Trapping Theorems. The proofs of these two theorems are omitted since, once this work done,
   both generalizations are straightforward. In  Section~$4$, we show that, after pulling back
    the metric to $\mathbb{R}^3$, the only possible
   leaves of the foliation are planes, cylinders and spheres. This forth section, involving
   a succession of simple had-hoc arguments, is crucial for its consequences although the
    proof itself may be skipped at first reading. Using the results exhibited
    in the two preceding sections, section~$5$ is
devoted to the proof of the $3$D modified
   conjecture. Finally a counter-example to the three dimensional general
   form of Turbiner's conjecture is exhibited in section~$6$.\\

\section{General setting}
In this section, we introduce the framework and the notions
necessary to prove the three dimensional version of the modified
Turbiner's conjecture. The two dimensional version, see
\cite{Imprimitively} and \cite{MilsonThesis} for a complete proof,
is also discussed.
  Recall that a Schrödinger operator on a $n$-dimensional
 Riemannian manifold is a second order
 differential operator  of the form
\[ \hh:=-\frac{1}{2}\Delta+U,\] where $\Delta$ is the
Laplace-Beltrami operator and $U$ is the potential function for the
physical system under consideration. If $g^{(ij)}$ is the
contravariant metric in a local coordinate chart and $g$ its
determinant,
 the operator is given by\[\hh=-\frac{1}{2}\sum_{i,j=1}^n
[g^{ij}\partial_{ij}+\partial_i(g^{ij})\partial_j-\frac{g^{ij}\partial_i(g)}{2g}\partial_j]+U.
\]

The main purpose of this paper is to study separation of variables
for operators that are equivalent to Schrödinger operators. To this
end we want to consider an adequate notion of equivalence, one that
preserves the formal spectral properties of the operators under
consideration. The appropriate notion, which will be used throughout
our work, is the following. Two differential operators are
\emph{locally equivalent} if there is a gauge transformation
$\hh\rightarrow~e^\lambda \hh e^{-\lambda}$ and a change of variable
that relate one to the other. With this notion of equivalence in
hand, one can in principle, verify if a general second order
differential operator $\hh_0$ is equivalent to a Schrödinger
operator. Indeed, every second order differential operator can be
given locally by
\[ \hh_0=-\frac{1}{2}\sum_{i,j=1}^n g^{ij}\partial_{ij}+\sum_i^n
h^i\partial_i+U.\] If the contravariant tensor $g^{(ij)}$ is
non-degenerate, that is if $g$ does not vanish, the operator can be
expressed as\begin{equation}\label{lap+vec+pot}
\hh_0=-\frac{1}{2}\Delta+\vec{V}+U,\end{equation} where $\vec{V}$ is
a vector field that has to be eliminated by a further
transformation, if possible.
 In this context, the equivalence
condition can be restated as follows: $\hh_0$ is gauge equivalent to
a Schrödinger operator if and only if $\vec{V}$ is  a gradient
vector field, with respect to the metric $g^{(ij)}$. For obvious
reasons, this criteria is called \emph{closure condition} and if
$\vec{V}=\grad(\lambda)$, the gauge factor has to be
$e^\frac{\lambda}{2}$. \\

 A differential operator $\hh_0$ is \emph{Lie algebraic} if it
is an element of the universal enveloping algebra of a finite
dimensional Lie algebra of first order differential operators. Given
a representation of a finite dimensional Lie algebra $\g$  by vector
fields on a manifold $\mathbf{M}$, one constructs a representation
by first order operators in the following way, see \cite{const}
 for details. For each element $a \in \g$, we define a first order differential operator
\[  T_a=a^\pi+\eta(a),
\] where $a^\pi$ is a vector field on $\mathbf{M}$ and $\eta$ is an element of $H^1(\g,
\mathcal{C}^\infty(\mathbf{M}))$. Thus for  $\{a_1,...a_m\}$, a
basis for $\g$,  a second order Lie algebraic differential operator
is given by
\begin{equation}\label{lie-alge} \mathcal{H}_0=\sum_{i,j=1}^mC^{ij}T_{a_i}T_{a_j}+\sum_{k=1}^mL^kT_{a_k}, \end{equation}
where $C^{ab}$ and $L^{c}$ are real numbers and without loss of
generality $C^{ab}=C^{ba}$. The class of Lie algebraic operators is
closed under gauge transformation, hence it make sense to consider
Lie algebraic operators that are equivalent to Schrödinger
operators. Again, if the induced contravariant metric $g^{(ij)}$ is
non-degenerate, the closure condition can be easily verify. This
condition can be written in term of the initial coefficients
$C^{ab}$ and $L^c$, however this
leads to complicated PDE's.\\

We now give an example to better understand these notions. Later on,
this same example will be used to illustrate the conjecture.
Consider $\mathbf{M}$ a three dimensional manifold diffeomorphic to
$\mathbb{R}^3$ and the Lie algebra
$\g=\mathfrak{a}_1\oplus\mathfrak{a}_1\oplus\mathfrak{a}_1$  spanned
by the first order differential operators
\[
\begin{array}{cccccc} T_1=\partial_u, & T_2=u\partial_u, &
T_3=\partial_v, & T_4=v\partial_v, & T_5=\partial_w, &
T_6=w\partial_w.
\end{array} \]
We define a Lie algebraic operator $\hh_0$  with the following
choice of coefficients: \[C^{ab}=\left(
  \begin{array}{cccccc}
    1/2 & 0 & 0 & 0& 0 & 0 \\
    0 &0 & 1 & 0 & 1 & 0 \\
    0 & 1 & 0 & 1 & 0 & 0 \\
    0 & 0 & 1 & 0 & 2 & 0 \\
    0 & 1 & 0 & 2 & 0 & 1 \\
    0 & 0 & 0 & 0 & 1 & 0 \\
  \end{array}
\right),\ \ \ L^c=\left(
         \begin{array}{c}
           0 \\
           2\alpha \\
           4\beta-4 \\
           4\alpha \\
           4\beta+4\gamma-6 \\
           4\alpha \\
         \end{array}
       \right),
\]
where $\alpha$, $\beta$ and $\gamma$ are real numbers. In term of
the $(u,v,w)$ coordinates, the operator reads as follows
\begin{eqnarray*} \hh_0&=&\frac{1}{2}\partial_{uu}
+2u\partial_{uv}+2u\partial_{uw}+2v\partial_{vv}+4v\partial_{vw}+2w\partial_{ww}\\
& &+2\alpha u\partial_u+(4\beta-3+4\alpha
v)\partial_v+(4\beta+4\gamma-5+4\alpha w)\partial_w.
\end{eqnarray*}
The metric associated to this operator is,
\begin{equation}\label{metricex1}
g^{(ij)}=-\left(
  \begin{array}{ccc}
    1 & 2u & 2u \\
    2u & 4v & 4v \\
    2u & 4v & 4w \\
  \end{array}
\right).
\end{equation}

If we forget the degeneracy issue for a moment, we can rewrite the
operator $\hh_0$ in term of the Laplace-Beltrami operator associated
to the metric (\ref{metricex1}). We obtain
\[
\hh_0=-\frac{1}{2}\Delta+2\alpha u\partial_u+(4\beta+4\alpha
v)\partial_v+(4\beta+4\gamma+4\alpha w)\partial_w.
\]
A direct calculation shows that the undesirable first order term can
be expressed as the gradient, with respect to (\ref{metricex1}), of
the scalar function \[\lambda=\alpha w+\beta \ln|v-u^2|+\gamma
\ln|w-v|.\] The closure condition is then satisfied and, scaling
with the factor $e^\frac{\lambda}{2}$, the operator constructed is
gauge equivalent to the following Schrödinger operator:
\[
\hh=-\frac{1}{2}\Delta+\alpha(2\beta+2\gamma+3)+\alpha^2w+\frac{\beta(\beta-1)}{v-u^2}+\frac{\gamma(\gamma-1)}{w-v}
.\]

 Observe that the
tensor $g^{(ij)}$ fails to be of full rank when
$g=16(w-v)(u^2-v)=0$. The inverse tensor $g_{(ij)}$ is singular on
the sets $\{w=v\}$ and $\{v=u^2\}$, thus the inner product is not
defined everywhere. We therefore have to allow degeneracy for the
contravariant metric. However, despite this flexibility, we want the
metric to behave reasonably well on the degeneracy locus. For this
reason, we introduce a generalization of the pseudo-Riemannian structure.\\

For $\mathbf{M}$ a real, analytic manifold and $g^{(ij)}$ a type
$(2,0)$ tensor field, we denote $\mathbf{D}_{g^{}}$ the locus of
degeneracy of the tensor. The analyticity requirement implies that
$\mathbf{D}_{g^{}}$ is either empty, a codimension $1$ subvariety or
$\mathbf{M}$. We set $\mathbf{M}_0=\mathbf{M} \backslash
\mathbf{D}_{g^{}}$ and we assume that $g^{(ij)}$ is not identically
degenerate. Thus $\mathbf{M}_0$ is an open, dense subset of
$\mathbf{M}$ and the connected components of $\mathbf{M}_0$ are
pseudo-Riemannian manifolds with boundary in $\mathbf{D}_{g^{}}$.\\

The pair $(\mathbf{M}, g^{(ij)})$ is called an
\emph{almost-Riemannian manifold} if for all pairs $u,v$ of analytic
vector fields with non-degenerate plane section $u\wedge v$ on
$\mathbf{M}_0$, the sectional curvature function $K(u\wedge v)$ has
removable singularities on $\mathbf{D}_{g^{}}$. Remark that if the
sectional curvature is constant on connected
components, which is the case if the Riemmanian curvature is zero, then $\mathbf{M}$ is an almost-Riemannian manifold.\\

Throughout this paper, we will focus on components of $\mathbf{M}_0$
for which the metric is positive definite and for which the
Riemannian curvature tensor is zero. For instance, in the previous
example, $g^{(ij)}$ is positive definite on
$\mathbf{R}=\{(u,v,w)|u^2<~v<w\}$  and
 we can easily check that the Riemannian curvature vanishes identically.\\

We will see that separation of variables is, in our context, closely
related to foliations by geodesics. Thus
 we would prefer to work in Euclidean geometry, where the geodesics are straight lines, instead of working on $\mathbf{M}$, a
 flat analytic manifold.
 The $3$D-Tiling Theorem will help us to achieve this by showing that, under certain conditions, there exists a global real-analytic map from $\mathbb{R}^3$ to $\mathbf{R}\subset \mathbf{M}$
where the contravariant metric of $\mathbf{R}$ is the pushforward of
the Euclidean metric. The following definitions and propositions
 will be necessary to establish this theorem.\\

To better understand the overall behaviour of the manifold around
the degenerate points, we need to quantify the degeneracy of the
contravariant metric. The degenerate points can be break up into two
categories. A point $p\in \mathbf{D}_g$ is called \emph{unreachable}
if all smooth curves with end points $p$
 have infinite length in the metric $g_{(ij)}$. Conversely a degenerate point is called \emph{reachable} if it can be attained by a finite length curve.
 If $\gamma(t):(0,1)\rightarrow \mathbf{M}_0$ is a geodesic segment, we denote $T$ the largest number, possibly $\infty$, such that
 $\gamma(t)$ can be extended by a geodesic with domain $(0,T)$. For $\mathbf{R}$ an open connected component of $\mathbf{M}_0$, we
 say that $\mathbf{M}$ is \emph{complete within $\mathbf{R}$} whenever for all geodesic segment lying within $\mathbf{R}$,
  either $T=\infty$, or $\lim_{t\rightarrow T}\gamma (t)$ is a reachable boundary point of $\mathbf{R}$. This extends the
  notion of completeness to almost-Riemannian manifold. The following Proposition will be useful.

 \begin{Pro}
Suppose the signature of $g^{(ij)}$ is positive definite within
$\mathbf{R}$, and that $\mathbf{R}$ is contained in a compact subset
of $\mathbf{M}$. Then $\mathbf{M}$ is complete within $\mathbf{R}$.
 \end{Pro}

 The degenerate
points are given by the zero set of the determinant of a $3$ by $3$
matrix which is, in general, not easy to handle. To circumvent this
issue, we will assume that either
\begin{enumerate}
\item $g^{(ij)}$ can be diagonalised,
\item $\mathbf{M}$ is a transverse, type changing manifold.
\end{enumerate}

For the  case $(1)$, there exists locally a coordinate system for
which $g^{(ij)}$ is expressed as
\begin{equation}\label{gij} g^{(ij)}=\left(
            \begin{array}{ccc}
              P(x,y,z) & 0 & 0 \\
              0 & Q(x,y,z) & 0 \\
              0 & 0 & R(x,y,z)\\
            \end{array}
          \right).\end{equation}
Thus, the determinant is given by the simple equation $g=PQR$ and
 we can assume without loss of generality that the metric is
degenerate at the origin. We define the \emph{order} of an analytic
function to be the smallest total degree of all the monomials with a
non-zero coefficient in its Taylor development. Thus,
 the order of $g$ will be the sum of the orders of the diagonal components of (\ref{gij}). Note that the smallest the order
 of $g$ is, the closest the metric is from being non-degenerate at the origin.\\

The next requirement will be needed to prove the three dimensional
version of the Trapping and Tiling Theorems when the metric is
diagonal. This condition does not appear in the two dimensional case
and we do not know yet if it is necessary. We say that a
contravariant metric tensor $g^{(ij)}$ given as (\ref{gij}) is
\emph{\avoir} if the components of the diagonal do not share
non-invertible factors. For instance the metric\[ g^{(ij)}=\left(
\begin{array}{ccc}
                        (1+x)^2 & 0 & 0 \\
                        0 & (1+x)y & 0 \\
                        0 & 0 & xz \\
                      \end{array}
                    \right)\] is \avoir \ while the metric
\[g^{(ij)}=\left(
  \begin{array}{ccc}
    x^2 & 0 & 0 \\
    0 & xy & 0 \\
    0 & 0 & xz \\
  \end{array}
\right)\] is not.\\

For the case $(2)$, it can be deduced that the metric can be
diagonalised and its diagonal form is \avoir\!. Indeed, recall that
$\mathbf{M}$ is a transverse, type changing analytic
 $m$-dimensional manifold if $\mathbf{M}$ is an
analytic manifold with a contravariant metric  $g^{(ij)}$ such that
at any point $x$ in the degenerate locus $\mathbf{D}_g$, we have:
\begin{enumerate}
\item $ d(det(g^{(ij)}))\mid_x\neq 0 $ for some (and hence any)
coordinate system,
\item $Rad_x:=\{ v_x\in T^*_x\mathbf{M}\mid g^{(ij)}(v , \cdot)=0 \}$ is transverse to
$T^*_x\mathbf{D}_g$,
\end{enumerate}

One can show, see \cite{tc} for details, that around any degenerate
point, there exists local natural coordinates
 $\{x^1,...,x^m\}$ such that
\[g^{(ij)}=
\left(\begin{array}{cc}
  g^{(ab)}
  & 0\\
  0& x^m
 \end{array}\right),\]
where $g^{(ab)}$ is non degenerate. In the three dimensional case,
$g^{(ab)}$ is a two by two matrix and can be diagonalised into
invertible functions. This leads us to a contravariant metric
\[g^{(ij)}=
\left(\begin{array}{ccc}
  P(x,y,z)
  & 0&0\\
  0& Q(x,y,z)&0\\
  0&  0& z
 \end{array}\right),\]
 for which the \avoirity property is satisfied.\\

Once the $3$D-Tiling Theorem have been established, there will be
one last major property needed in our study: the imprimitivity of
the action. Since we are dealing with a Lie algebraic operator, we
can
 assume that the domain of the operator is an homogeneous space $\mathbf{M}=\mathbf{G}/\mathbf{H}$ where $\g$
 is the Lie algebra corresponding to $\mathbf{G}$. Recall that the action of $\mathbf{G}$ on $\mathbf{M}$ is \emph{imprimitive}
  if there exists a foliation of $\mathbf{M}$ that
 is invariant under the action of $\mathbf{G}$.
 In three dimensional Euclidean space, the invariant leaves can be either curves or surfaces,
see \cite{Lie} for a more detailed description of the possible
leaves. Throughout this paper we will only consider foliations by
surfaces, this type of action will be called \emph{\imprimitive
action}. Note however that it would be
 very interesting to study the case of a foliation by curves.
  In infinitesimal terms, if the leaves of the foliation are given by $\{ \lambda=\textrm{constant} \}$, then
 $a^\pi(\lambda)=f(\lambda)$ for every element of the lie algebra $\g$. This second criteria can be generalized
 to extend the notion of \imprimitivity to differential operators.  If the level sets of the function $\lambda$ are the leaves of the
 foliation, the operator $ T_\alpha$ is said to acts \imprimitively  if $ T_\alpha (\lambda)$ and $\lambda$ are
 functionally
  dependent. One can easily show, see \cite{Imprimitively} for details, the following
  \begin{Pro}\label{imp} If the operators $\{ T_a  :\ a \in g\}$ act \imprimitively, then there is a $\mathbf{G}$-invariant
   foliation on $\mathbf{M}$.
  \end{Pro}

The central point here is that Lie algebraic operators generated
 by these \imprimitive generators will behave the same way.
Indeed the operator $\hh_0$ applied to $\lambda$ will give back a
function of $\lambda$. By taking $\lambda$ as coordinate, the
operator
 will separate in that variable and we will show that the equivalent Schrödinger
 operator $\hh$ will also separate partially.\\

One of the key arguments for the final theorem is that the invariant
foliation $\Lambda$ is perpendicular to a geodesic foliation. This
is due to the Lie algebraic construction of the metric $g^{(ij)}$,
and, according to the 3D-Tiling Theorem, this perpendicular
foliation can be pull back to the Euclidean space where the
geodesics are well known: straight lines. We will not go into the
details, everything being already exhibit in \cite{Imprimitively}
and \cite{MilsonThesis}, but we will state the mains results
necessary to prove the theorem.\\

We will denote $\Lambda^\perp$ the distribution of tangent vectors
that are perpendicular to $\Lambda$. For a Lie algebraic metric with
invariant foliation $\Lambda$, one can prove the following:

\begin{Theo}If $\Lambda^\perp$ is tangent to a geodesic of $M$
 at one point, then the geodesic is an integral manifold of $\Lambda^\perp$.
\end{Theo}

In the context of the modified three dimensional Turbiner conjecture
$\Lambda$ is a rank $2$ $\mathbf{G}$-invariant distribution, thus,
being a rank $1$ distribution, $\Lambda^\perp$ is necessarily
integrable. We then get:

\begin{Co}\label{perpendicular} If rank$(\Lambda^\perp)=1$, then the integral curves of
$\Lambda^\perp$ are geodesic trajectories.
\end{Co}

After an investigation of the possible foliations of $\mathbb{R}^3$
which are in accordance to our problem,
 we will be able to show that the partial separation will occurs either in Cartesian, cylindrical or spherical coordinates.
   Note that, as for the two dimensional case,
 the extra hypothesis are necessary. Based on a primitive action, an explicit flat Lie algebraic Schrödinger operator, for which there
 is no separation of variable, will be exhibited at the end of this paper.\\

 Let us now illustrate how the \imprimitivity of the action affects the operator constructed previously.
 Recall that if we consider the domain $\mathbf{R}\subset\mathbf{\mm_0}$ where the metric is positive definite, the operator
  reads as follows
\[\hh_0=\Delta+\grad(\alpha w+\beta
\ln|v-u^2|+\gamma \ln|w-v|).
 \]
This operator, and its equivalent Schrödinger operator, illustrate
clearly that separation arises from invariant foliations,
 indeed the separation occurs in each of
the 3 possible systems of coordinates. We shall not expect this in
general. The three separations reflect the fact that the group
action allows not only one but three distinct invariant foliations:
\[ \{ u=\textrm{const.} \},\ \ \ \ \    \{  v=\textrm{const.}  \}, \ \ \ \ \  \{  w=\textrm{const.} \}.
 \] It is now guaranteed that  $\hh_0(\lambda)=f(\lambda)$ for $\lambda \in \{u,v,w \}$.\\

In term of Cartesian coordinates $(x,y,z)$, the original coordinates
are given by
\[ u=x,  \ \ \ \ v=x^2+y^2 \ \ \ \ w=x^2+y^2+z^2.
\]
Thus the leaves of the foliation are planes, cylinders and spheres.
For each of these foliations we will consider respectively the
Cartesian, cylindrical $(r,\theta, z)$, and spherical $(r,\theta,
\phi)$ coordinates.
 Hence, in the Cartesian system, the coordinate $x$ separate in the operator, and for the two other systems
 the radial coordinate can be
separated.   An other extra property of this
  operator is the fact that, for each of these coordinate systems, the operator also separate in the two other coordinates.\\

In these three coordinate systems, the operator is given by
\begin{eqnarray*}\hh_0&=&\Delta+\grad(\alpha (x^2+y^2+z^2)+\beta
\ln|y|+\gamma \ln|z|),
 \\
\hh_0&=&\Delta+\grad(\alpha (r^2+z^2)+\beta
(\ln|r|+\ln |\sin\theta|)+\gamma \ln|z|),\\
\hh_0&=&\Delta+\grad(\alpha (r^2)+\beta
(\ln|r|+\ln|\sin\phi|+\ln|\sin \theta|)+\gamma( \ln|r|+\ln|\cos
\theta|)).
 \end{eqnarray*}
By applying the operator $\hh_0$ to
$\Phi(x_1,x_2,x_3)=\Phi_1(x_1)\Phi_2(x_2)\Phi_3(x_3)$, one easily
verifies
 that the equation separate in three equations, each of them involving only one variable. After the required gauge transformation, the Schrödinger operator reads as
\begin{eqnarray*}
\hh&=&\Delta+\alpha (2\beta+2\gamma+3)+\alpha^2(x^2+y^2+z^2)+\frac{\beta(\beta-1)}{y^{2}}+\frac{\gamma(\gamma-1)}{z^{2}},\\
\hh&=&\Delta+\alpha
(2\beta+2\gamma+3)+\alpha^2(r^2+z^2)+\frac{1}{r^2}[\frac{\beta(\beta-1)}
 {\sin^2\theta}]+\frac{\gamma(\gamma-1)}{z^{2}},\\
 \hh&=& \Delta+\alpha (2\beta+2\gamma+3)+\alpha^2(r^2)+\frac{1}{r^2\sin^2\theta}[\frac{\beta(\beta-1)}
 {\sin^2 \phi}]+\frac{1}{r^2}[\frac{\gamma(\gamma-1)}{cos^2\theta}].
\end{eqnarray*}Once again, the three operators separate in their respective coordinate systems since
the three potentials satisfy the separation condition, see \cite{sep} for more details.\\

The aim of the next sections is to prove the three dimensional
version of the following modified
Turbiner's conjecture proved by Rob Milson in \cite{Imprimitively}:\\

\begin{Theo}
Let $\hh_0$ be a second-order Lie algebraic operator generated by
the $T_a$ as per $(\ref{lie-alge})$, $g^{(ij)}$ the induced
contravariant metric and $\mathbf{R}$ a connected component of
$\mathbf{\mm_0}$ for which $g^{(ij)}$ is positive definite.
 Suppose the following statements are true:
\begin{enumerate}
\item $\hh_0$ is gauge equivalent to a Schrödinger operator;
\item $(\mathbf{R}, g^{(ij)})$ is isometric to a subset of the Euclidean plane;
\item The operators $\{ T_a : a \in \g \}$ act imprimitively;
\item $\mathbf{R}$ is either compact, or can be compactified in such a way that the $\mathbf{G}$-action on
$\mathbf{R}$ extends to a real-analytic action on the
compactification.
\end{enumerate}
Then, both the eigenvalue equation $\hh_0 \psi=E\psi$, and the
corresponding Schrödinger
equation separate in either a Cartesian, or a polar coordinate system.\\
\end{Theo}
To this end, we will follow a path which is similar to the one
 followed by Rob Milson. However, as mentioned
above, the extra requirement that the metric is diagonalisable and
\avoir \!\!, will be required.

\section{Trapping and Tiling} The main objective of this section is
to show that, under some conditions, there is a global map from the
Euclidean space to the positive definite region of a flat three
dimensional almost-Riemannian compact manifold. As for the planar
case, the $3$D-Tiling Theorem will follow principally from the
$3$D-Trapping Theorem. That later assure that the flow of a gradient
vector field can never cross the locus of degeneracy. Note that
through this section, the contravariant metric will be taken to be
diagonal and the \avoirity
property will be needed to prove both theorems.\\
\subsection{$3$D-Trapping theorem} The trapping property is a feature
shared by every flat diagonal \avoir almost-Riemannian metric whose
coefficients are analytic functions. All the work involved in the
proof is based on an appropriate expression of the diagonal
components of the Ricci curvature tensor. We will use the local
coordinates $(x,y,z)$, that we will sometimes denote $(x^1,x^2,x^3)$
to ease the notation. Thus, we define
$H^i=\sum_{j}g^{ij}\frac{\partial}{\partial
x^i}=g^{ii}\frac{\partial}{\partial x^i}$ and evaluate the diagonal
components of the Ricci curvature tensor using
 the frame $\{ H^1, H^2, H^3 \}.$   After some work of
simplifications and rearrangements, we obtain the
 following three expressions:
 \begin{eqnarray*}
 2(R_{11})g^2&=&-3(H^1(g))^2+2g(H^1(H^1(g))\\
 &&+g^2 [P_yQ_y+P_zR_z+2QP_{yy}
 +2RP_{zz}+P_x^2-2PP_{xx}]\\
&&+g [2P^3Q_xR_x+P^2QP_xR_x+P^2RP_x Q_x-
 PQ^2P_yR_y\\
 &&-PR^2P_zQ_z
 -3QR^2P_z^2-3Q^2RP_y^2],\\
 &&\\
 2(R_{22})g^2&=&-3(H^2(g))^2+2g(H^2(H^2(g))\\
 &&+g^2 [P_xQ_x+Q_zR_z+2RQ_{zz}
 +2PQ_{xx}+Q_y^2-2QQ_{yy}]\\
 &&+g [2Q^3P_yR_y+PQ^2Q_yR_y+Q^2RP_yQ_y
-P^2QR_xQ_x\\
&&-QR^2P_zQ_z
 -3PR^2Q_z^2-3P^2RQ_x^2],\\
 &&\\
2(R_{33})g^2&=&-3(H^3(g))^2+2g(H^3(H^3(g))\\
 &&+g^2 [Q_yR_y+P_xR_x
 +2QR_{yy}+2PR_{xx}+R_z^2-2RR_{zz}]\\
&&+g [2R^3P_zQ_z+QR^2P_zR_z+PR^2Q_z R_z-
 Q^2RP_yR_y\\
 &&-P^2RQ_xR_x
 -3PQ^2R_y^2-3P^2QR_x^2].
 \end{eqnarray*}

\begin{Pro}\label{pretrapping} Let $g^{ij}$ be a diagonal, \avoir three dimensional contravariant metric tensor
 with analytic coefficients. If the diagonal
elements of the Ricci curvature tensor are identically zero, then
there exists locally defined, analytic functions $\mu^1$, $\mu^2$
and $\mu^3$ such that
\[  H^i(g)=\mu^i\cdot g \  \textrm{ for } i=1,2,3.\]
\end{Pro}

\noindent{\sc {\bf Proof: }} Obviously, such functions exist around
points where the determinant does not vanish. We can assume that $g$
is zero at the origin and we focus on $H^1$ first. The ring of
convergent power series with complex coefficients is a unique
factorization domain, thus, up to multiplication by invertible
functions, $g$ factors uniquely into a product of irreducible,
complex valued, analytic functions that are zero at the origin.  Let
$f$ be such factor, and let $k$ be its multiplicity, i.e.
$g=f^k\sigma,$  with $f$ and $\sigma$ coprime. Since $g^{(ij)}$ is
\avoir \!\!, $f^k$ is only a factor of one of the diagonal elements
and if $k$ is greater then one, $f^{k-1}$ divides  the three partial
derivatives of this component. Thus, one easily sees that $f^{2k-1}$
is a factor of the two last summands of $2(R_{11})g^2$:
\[ g^2[P_yQ_y+...-2PP_{xx}]+ g[2P^3Q_xR_x+...-3Q^2RP_y^2]. \]
Since $R_{11}$ is
 identically zero,  the remaining summand,
 $3(H^1(g))^2-2g(H^1(H^1(g))$,
 must also be divisible by $f^{2k-1}$. But, the preceding term can be written as
 \[k(k+2)\sigma^2(H^1(f))^2f^{2k-2}+\rho f^{2k-1}\]
where $\rho$ is some analytic function. Thus
$k(k+2)\sigma^2(H^1(f))^2$ must be divisible by $f$. Recall that
$\sigma$ is relatively prime to $f$ and $k(k+2)>0$, which force
$H^1(f)$ to be divisible by $f$. The same must be true for all
non-invertible irreducible factors of $g$, (and obviously true for
the invertible factors), therefore $H^1(g)$ is divisible by $g$. The
same argument
holds for $H^2$ and $H^3$.~\hspace{\stretch{1}} $\Box$\\

Note that, without the \avoirity requirement, the first term of the
last summand of $2(R_{11})g^2$ is only guaranteed to be divisible by
$f^{2k-2}$ which does not allows us to establish the claim. However,
maybe an other rearrangement of the terms could lead to the same
conclusion without this extra hypothesis.\\

 From this proposition,
the 3D-Trapping Theorem follows immediately. Being identical to the
one given in Corollary $6.4.2$ of \cite{MilsonThesis}, the proof is
omitted.

\begin{Theo}\label{trapping}(The 3D-Trapping Theorem) Let $g^{ij}$ be as in the preceding
theorem, and let $f$ be an analytic function. Then the flow of
$\grad(f)$ can never cross the locus of degeneracy. More precisely,
the trajectories of the flow of $\grad(f)$ are either contained in
the locus of degeneracy of $g^{ij}$, or never intersect it.
\end{Theo}

\subsection{$3$D-Tiling theorem}

In what follows, using the 3D-Trapping Theorem, we will prove a
three dimensional version of Rob Milson's Tiling Theorem. As before,
$\mathbf{M}$ is a compact, three dimensional, almost-Riemannian
manifold endowed with $g^{(ij)}$ a \avoir and flat metric with
diagonal analytic coefficients.  $\mathbf{R}$ is a
region where the metric is positive definite. \\

The key argument for this proof is that either the degenerate points
are unreachable, either the metric $g^{(ij)}$ is the push-forward of
a non-degenerate metric $\tilde{g}^{(ij)}$. But before proving this
proposition, the two following lemmas, deduced from the Proposition
\ref{pretrapping} and the \avoirity property of the metric, will
simplify the subsequent work. Under the same hypothesis, we have the
following:

\begin{Le}\label{lemme1} If $f$ is a non-invertible,
irreducible factor of $g^{ii}$, then, for $i\neq j$, $f$ is a factor
of $g^{ii}_{x^j}$ and a factor of $f_{x^j}$.
\end{Le}
\noindent{\sc {\bf Proof: }} Suppose that $P=f^k\sigma$, where $f$
is a non-invertible irreducible factor and $(f,\sigma)=1$. By the
\avoirity property of the metric, $f$ is also coprime to $Q$ and
$R$, and from Proposition~\ref{pretrapping},
\[ H^2(g)=Q(P_yQR+PQ_yR+PQR_y)=\mu^2\cdot PQR.\]
$P$ being a factor of all but one summands of the middle term,
$P_yQ^2R$ has to be also divisible by $f^k$, forcing $f^k$
 to divide $P_y$. Furthermore,
\[ P_y= \left\{  \begin{array}{ll}kf^{k-1}f_y \sigma+f^k \sigma_y& \textrm{ if } k>1,\\
f_y& \textrm{ if } k=1,
\end{array}\right. \]
 thus $f$ needs to be a factor of $f_y$. \hspace{\stretch{1}} $\Box$\\

\begin{Le} \label{lemme2} Given $g^{ii}$, a
 diagonal component of the contravariant
metric $g^{ij}$, its non-invertible factors are functions of the
variable $x^i$ only.
\end{Le}

\noindent{\sc {\bf Proof: }}Consider $f$, a non-invertible factor of
$R$. From the analycity requirement, $f$ can be expressed locally by
the following convergent power series,
\[ f=\sum_{i,j,k=0}^{\infty}f_{ijk}x^iy^jz^k, \ \textrm{ where } \ f_{000}=0.\]
According to  Lemma~\ref{lemme1}, $f_x=f\cdot h$, for $h$ an
analytic function. The Taylor series of $f_x$ can therefore be given
as a product of two series:
\begin{equation}\label{series} f_x=\sum_{i,j,k=0}^{\infty}
if_{ijk}x^{i-1}y^jz^k=\sum_{i,j,k=0}^{\infty} f_{ijk}x^iy^jz^k \cdot
\sum _{a,b,c=0}^{\infty}h_{abc}x^ay^bz^c.\end{equation} If we
suppose that there exist positive integers $i$, such that $f_{ijk}
\neq 0$, we can
 fix $(\alpha, \beta, \gamma)$,
  the smallest triple
(with respect to the lexicographic order) such that $f_{\alpha \beta
\gamma}\neq 0$.  The coefficient of the monomial
$x^{\alpha-1}y^{\beta} z^{\gamma}$ is $\alpha f_{\alpha \beta
\gamma}$, and according to (\ref{series}), it can also be given by
\[\sum_{\tiny{\begin{array}{c}
i+a=\alpha-1\\j+b=\beta\\k+c=\gamma\end{array}}}f_{ijk}\cdot
h_{abc}.\] But all the coefficients $f_{ijk}$ are zero since
$i=\alpha-1-a<\alpha$. Consequently $f_{\alpha \beta \gamma}= 0$, a
contradiction. So, we have $f_{ijk}=0$ for all $i\neq 0$, and the
same argument is used to show that $f_{ijk}=0$ for all $j\neq 0$.
Therefore
\[ f=\sum_{k=0}^\infty f_{00k}z^k=f(z).\] \hspace{\stretch{1}} $\Box$\\

We can now prove the following strong criteria for the
unreachability of a degenerate point. As for the rest of this paper,
the degenerate point will be taken to be the origin.
\begin{Pro} If the order of one of the diagonal components $g^{ii}$ is greater than one,
then the origin is unreachable.
\end{Pro}
\noindent{\sc {\bf Proof: }} Without loss of generality,
$R=z^l(1+f(x,y,z))$ where $l>1$. We will compare the metric
$g^{(ij)}$ to
\[ \tilde{g}^{(ij)} =\left(
      \begin{array}{ccc}
        1 & 0 & 0 \\
       0 & 1 & 0 \\
        0 & 0 & z^2\\
      \end{array}
    \right),
\] which is a flat metric whose origin is known to be unreachable. We can write the
contravariant metric as
\[ g^{(ij)} =\left(
      \begin{array}{ccc}
        P & 0 & 0 \\
       0 & Q & 0 \\
        0 & 0 & z^2\tilde{R}\\
      \end{array}
    \right), \]
    where $P,Q,\tilde{R}$ are non-singular at the origin. We can
    find a neighborhood $N$ and a upper bound $K>0$ such that
    $\sup_N \{P,Q,\tilde{R} \}\leq K$. If we consider the region $\mathbf{R}\cap
    N$, we must have
    \[ \langle v,v \rangle_g \geq \frac{1}{K}\langle v,v\rangle_{\tilde{g}},\]
    for all tangent vectors $v$. Indeed
    \[
    \langle v,v \rangle_g=\frac{v_1^2}{P}+\frac{v_2^2}{Q}+\frac{v_3^2}{z^2\tilde{R}}\geq
    \frac{v_1^2}{K}+\frac{v_2^2}{K}+\frac{v_3^2}{z^2K}=\frac{1}{K}\langle v,v\rangle_{\tilde{g}}.\]
The length functional on curves in the metric $g$ is bounded below
by $\frac{1}{K}$ time the length functional in the metric
$\tilde{g}$. The origin, unreachable with respect to $\tilde{g}$, is
therefore
 unreachable with respect to
$g$ as well.~\hspace{\stretch{1}} $\Box$\\

\begin{Co} If the origin is reachable, then the order of $P,Q$ and $R$ is
at most one.
\end{Co}
This leads us, up to relabeling of the variables, to three
possibilities:
\begin{eqnarray}
g=PQR&=&z(1+f(x,y,z))\label{cas1},\\
g=PQR&=&yz(1+f(x,y,z))\label{cas2},\\
g=PQR&=&xyz(1+f(x,y,z))\label{cas3},
\end{eqnarray}
that enable us to prove the following key lemma.

\begin{Pro}\label{pretiling} A degenerate point is either an unreachable point, or
there exists a contravariant, non-degenerate metric tensor
$\tilde{g}^{(ij)}$ with analytic coefficients defined on some
neighborhood $N\subset \mathbb{R}^3$ and an analytic map $\phi:
N\rightarrow R$ such that $\phi_*(\tilde{g})=g.$
\end{Pro}
\noindent{\sc {\bf Proof: }}If the origin is reachable, we are in
one of three previous possibilities, say the case (\ref{cas3}).
Since each diagonal component has order $1$, from
Lemma~\ref{lemme2}, we can write $P=4x\tilde{P}$, $Q=4y\tilde{Q}$
and $R=4z\tilde{R}$ where $\tilde{P},\tilde{Q},\tilde{R}$ are
invertible. We consider the analytic map
 given by
\[
\phi_3: = \left\{ \begin{array}{lll}
x & =& \xi^2 \\
y & =& \eta^2 \\
z & =& \mu^2
\end{array} \right. ,
\] and we take $N$, the domain of this map, to be a neighborhood of
the origin sufficiently small so that the image of the map is
contained in $R$. One easily verify that, via this map, $g^{(ij)}$
is the pushforward of

\[ \tilde{g}^{(ij)} =\left(
      \begin{array}{ccc}
        \tilde{P} & 0 & 0 \\
       0 & \tilde{Q} & 0 \\
        0 & 0 & \tilde{R}\\
      \end{array}
    \right), \] which is non-degenerate at the origin. The cases
    (\ref{cas1}) and (\ref{cas2}) are resolved the same way, by considering
    respectively the maps
    \begin{equation}
\phi_1: = \left\{ \begin{array}{lll}
x & =& \xi \\
y & =& \eta \\
z & =& \mu^2
\end{array} \right.  \textrm{and  }  \phi_2: = \left\{ \begin{array}{lll}
x & =& \xi \\
y & =& \eta^2 \\
z & =& \mu^2  \end{array} \right. .\end{equation} ~\hspace{\stretch{1}} $\Box$\\

The maps $\phi_i$ will be called \emph{$2^i$th-fold maps}. The name
reflects the fact that the $(\xi, \eta, \mu)$ space generically
covers the $(x,y,z)$ space in a $2^i$-to-one relationship. The
exception being the folding planes, $z=0$ for $\phi_1$, $y=z=0$ for
$\phi_2$ and $x=y=z=0$ for $\phi_3$. With this key lemma in hand, we
can now assert that the positive-definite region of the
almost-Riemannian manifold is isometric to the Euclidean space
modulo a discrete group isometric symmetries. The proof is identical
to the one given in \cite{MilsonThesis} for the two dimensional case
and is based on the fact that reachable degenerate points are, in a
way, removable. The three dimensional case is a little simpler since
the only analytic maps we need to consider are the $2^i$th-fold
maps. For these reasons, we will omit the proof.

\begin{Theo}(The 3D-Tiling Theorem)\label{Tiling} Let $M$ be a compact
three dimensional flat almost-Riemannian manifold with diagonal
\avoir metric. Then, there exists a globally defined, real-analytic
map $ \psi : \mathbb{R}^3 \rightarrow M$ such that $g^{(ij)}$ is the
push forward of the Euclidean metric, and such that $\psi$ covers
all of $\mathbf{R}$ plus the reachable portions of its boundary.
Furthermore, $\mathbf{R}$ is isometric to the quotient
 $\mathbb{R}^3 / \Gamma$, where $\Gamma$ is the group of isometries
 $\gamma$ such that $\psi=\psi \gamma$.\end{Theo}
Note that since $\psi$ is a $2^i$th-fold map, the group of
isometries
 is indeed the group of reflections along the folding planes.

\section{Foliations}
In this section, we intend to determine what are the possible rank
two foliations of $\mathbb{R}^3$ that are perpendicular to straight
lines. This intermediate result, used in conjunction with the
$3$D-Tiling Theorem, will be used to prove that the function
$\lambda$, whose level sets are the leaves of the invariant
foliation, is a coordinate of
either the Cartesian, the cylindrical or the spherical coordinates systems.\\

 The rank two leaves are complete on $\mathbf{M}$ but they may
cross the reachable part of the degenerate locus
$\mathbf{D}_{g^{}}$. Is is not clear a priori that the pull back of
these leaves is also complete in $\mathbb{R}^3$. Indeed, the rank of
these leaves may drop where the Jacobian of
$\psi:\mathbb{R}^3\rightarrow \mathbf{R}\subset \mathbf{M}$ is
degenerate. To avoid confusion, we denote $\mathbf{S}_{g^{}}$ the
degeneracy locus of the foliation in
$\mathbb{R}^3$ and we have the inclusion $\mathbf{S}_{g^{}}\subseteq\psi^{-1}(\mathbf{D}_{g^{}})$.\\

From the $3$D-Tiling Theorem, $\mathbf{R}$, the positive definite
region of the manifold,  is isometric to the quotient
 $\mathbb{R}^3 / \Gamma$ where $\Gamma$ is a discrete group of reflections.
 Thus $\mathbb{R}^3$ is tiled into isometric regions where the pull back
 of each leaf is repeated. If a rank two leaf of $\mathbf{R}$ crosses
 $\mathbf{D}_{g^{}}$, its pull back will be reflected on the other side of $\psi^{-1}(\mathbf{D}_{g^{}})$.
 Hence the rank two leaves in $\mathbb{R}^3$ can be extended without restriction but they
 might fail to be smooth. \\

 However, according to the $3$D-Trapping Theorem, the trajectories of
the flow of the gradient of
 $\lambda$ are either contained in the locus of degeneracy, or never intersect it. This
 forces the leaves to cross $\psi^{-1}(\mathbf{D}_{g^{}})$ perpendicularly, thus, we can conclude that
 the rank two leaves are also smooth in $\mathbb{R}^3$.\\

Therefore throughout this section,  $\Lambda$  will denote a
foliation of $\mathbb{R}^3$ which is of rank two almost everywhere.
By degenerate points we refer to $\mathbf{S}_{g^{}}$, the points
where the rank drops. One easily sees that the rank two leaves never
cross the locus of degeneracy. In accordance with Corollary
\ref{perpendicular}, the leaves of $\Lambda^{\perp}$, the
perpendicular foliation, are straight lines at every non-degenerate
point. Our aim is to show that the leaves of $\Lambda$ can only be
planes, infinite cylinders or spheres. Before proving this result we
need to establish some
notations together with two lemmas.\\

 For any point $x\in \mathbb{R}^3$, we
denote $\mathcal{M}_x$ its leaf and, for any curve $c$ contained in
a rank two leaf $\mm$, we denote $S_{c}$ the ruled surface generated
by the normals of $\mm$ along $c$. Throughout this section the
non-degenerate points will be dense and we will use the definitions
found in \cite{SPIVAK3} to describe solids.

\begin{Le}\label{strates}
Let $c(t)$ be a continuous family of curves parametrized by $t\in (
-\delta, \delta )$, contained in $\mathcal{M}$, a rank $2$ leaf.
Suppose that for every $t_1\neq t_2\in( -\delta, \delta )$, the
curves $c(t_1)$ and $c(t_2)$ are distinct almost everywhere.
 If all the surfaces
$S_{c(t)}$ intersect each others, then they all intersect at $\ii$
which is of dimension is at most one.
\end{Le}
\centerline{
 \epsfxsize=2,2in \epsfysize=2in \epsfbox{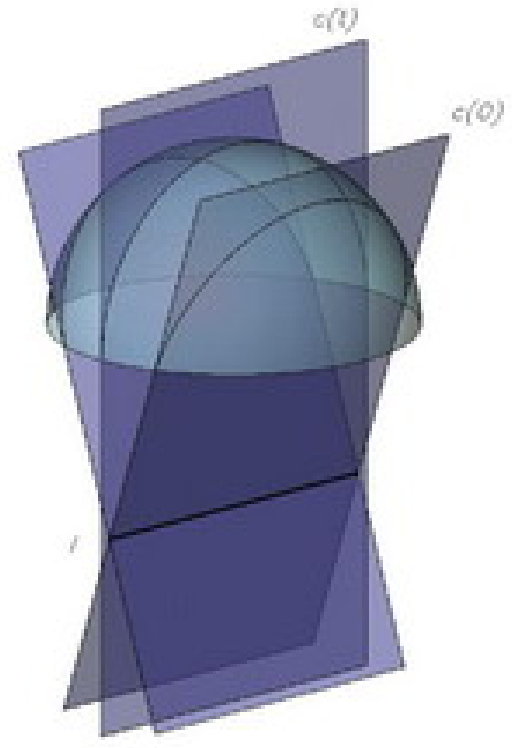}}
\noindent{\sc {\bf Proof: }}  Two different curves, $c(t_1)$ and
$c(t_2)$, can only intersect at points, hence the two ruled
surfaces, $S_{c(t_1)}$ and $S_{c(t_2)}$,  can intersect at most in a
one dimensional set. For any $\tau \in ( -\delta, \delta )$, we
define $\mathcal{I}_\tau$ the intersection of $S_{c(\tau)}$ with all
the other surfaces.
\begin{displaymath} \mathcal{I}_\tau:=\bigcup_{\scriptsize{\begin{array}{c}
               t\in ( -\delta,\delta ), \\
               t \neq \tau
\end{array}}}S_{c(t)}\cap S_{c(\tau)}\neq\emptyset.
\end{displaymath}
The set $\mathcal{I}_0$ can not be a surface. Otherwise, by
smoothness of the leaf,  $\mathcal{I}_\rho$ would also be a surface
for $|\rho |< \epsilon $ and $\epsilon$ sufficiently small. Thus we
would have
\begin{displaymath} \mathcal{I}:=\bigcup_{|\rho |<
\epsilon}\mathcal{I}_\rho,
\end{displaymath}
a three dimensional degenerate set. Hence, every surface $S_{c(t)}$
will intersect $S_{c(0)}$ at $\mathcal{I}_0$ which is of dimension
is at most one. The same argument holds for each set $\ii_t$.
Therefore, $\mathcal{I}_0=\mathcal{I}_t$ for all $ t \in ( -\delta,
\delta )$. \hspace{\stretch{1}} $\Box$\\

\begin{Le}\label{etoile}
Let $c(t)$ be a continuous family of curves parametrized by $t\in (
-\delta, \delta )$ and contained in $\mathcal{M}$, a rank $2$ leaf.
 If,  for all $t\in ( -\delta, \delta )$, each normal of $\mm$ along $c(t)$
 intersects a degenerate curve $\mathcal{I}$, then $\mathcal{I}$  is
parallel to $c(t)$ for all $t\in ( -\delta, \delta )$.\\
\end{Le}

\noindent{\sc {\bf Proof: }} Let $\mathcal{U}\subset \mathbb{R}^3$
be an open set contained in \begin{displaymath}S:=\bigcup_{t\in (
-\delta, \delta )}S_{c(t)}.\end{displaymath} We pick $x \in
\mathcal{U}$, a $\Lambda$-rank $2$ point, and, by completeness,
$\mathcal{M}_x$ is crossed perpendicularly by each normal associated
to the family $c(t)$. This section of the leaf, given by
$\bigcup_{t\in ( -\delta, \delta )}\mm_x \cap S_{c(t)}$, is
therefore parallel to $\mathcal{M}$. We denote $w(t)$ the
intersection of $\mathcal{M}_x$ with $S_{c(t)}$ and we note that
$c(t)$ is at constant distance from $w(t)$.

\centerline{
 \epsfxsize=1,4in \epsfysize=1,4in \epsfbox{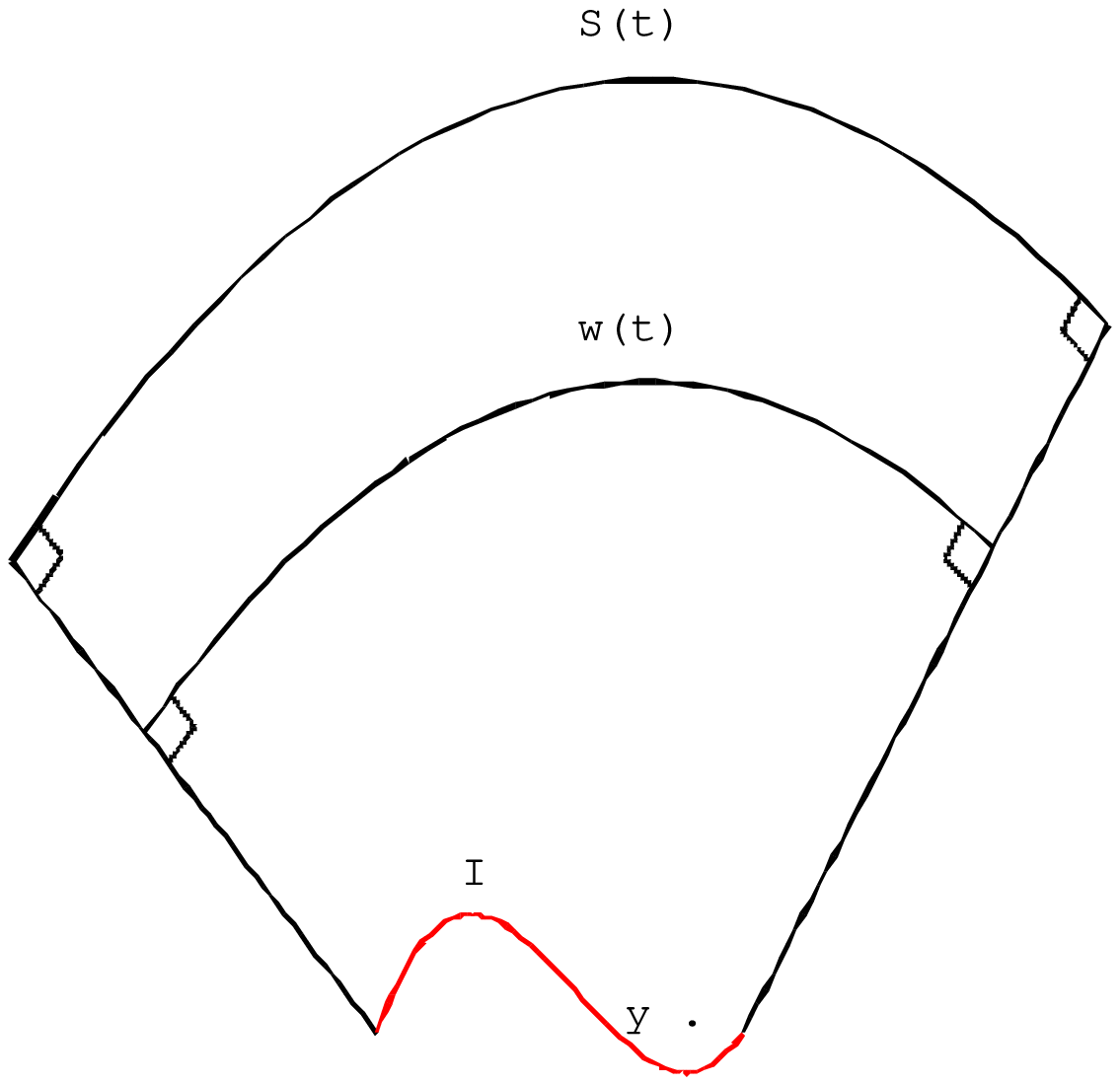}}
\noindent If  $\mathcal{I}$ was not parallel to a curve $c(t)$, one
could easily choose a rank $2$ point, say $y$, sufficiently closed
to $\ii$, for which $\mathcal{M}_y$ would intersects the degenerate
set $\mathcal{I}$. This is impossible since the leaf $\mathcal{M}_y$
has $\Lambda$-rank $2$
everywhere. \hspace{\stretch{1}} $\Box$\\

Remark that the intersection curve $\ii$, being parallel to the
leaves, has to be non-singular. Note also that the curves in the
family are necessarily all parallel to each others. We can now prove
the following three propositions, which, put together, will enable
us to conclude about the three possible foliations.

\begin{Pro}\label{sphere}Let $\mm$ be a rank $2$ leaf of the foliation $\Lambda$,
if there exists an open $\uu \subset \mm$ for which the Gaussian
curvature is positive, then $\mm$ is a sphere.
\end{Pro}
\noindent{\sc {\bf Proof: }} Let $y\in \uu$ and consider
$c(0)\subset \uu $ the segment of the curve starting at $y$ and
following the leaf $\mm$ in a given direction $\pm \vec{v}$. If we
fix the end points $a$ and $b$ and slide $c(0)$ in the two
directions perpendicular to $\vec{v}$, we get a family of curves
$c(t)$, $t\in ( -\delta, \delta )$, contained in $\uu$. We denote
$\vv$ the subset generated by the curves $c(t)$.

\centerline{
 \epsfxsize=1,9in \epsfysize1,5in \epsfbox{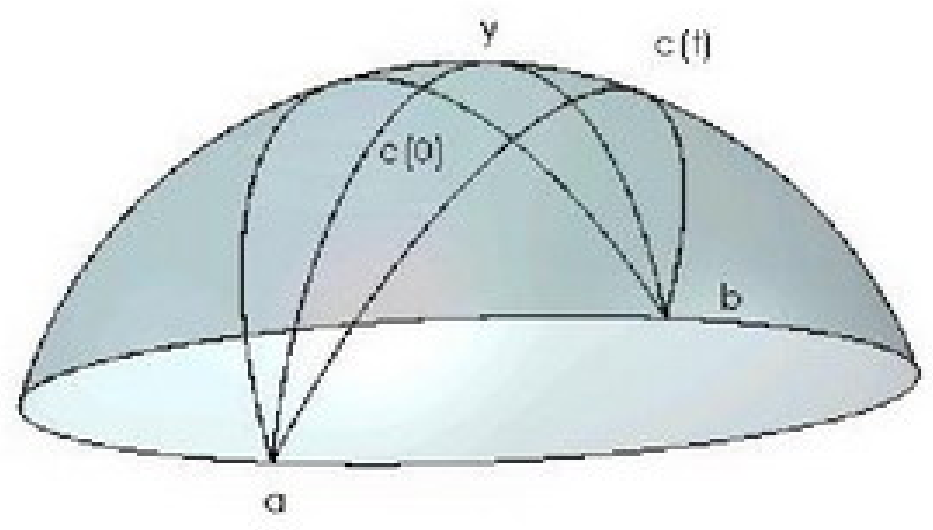}}

 Since the two principal curvatures are non-zero in $\uu$,
the normal surfaces $S_{c(t)}$ intersect in a connected component
and, by Lemma~\ref{strates} the intersection is either a point $q$,
either a connected curve $\ii$, parallel to $c(t)$ by
Lemma~\ref{etoile}. Suppose first that the intersection is a curve.
Being parallel to $\ii$, the surface
  $\vv$ has to be
contained in a twisted cylinder centered at $\ii$. For each point
$p$ on the curve $\ii$, we denote $c'(p)$, the intersection of $\vv$
with the normal plane of $\ii$ at $p$.

\centerline{
 \epsfxsize=2,5in \epsfysize=1,5in \epsfbox{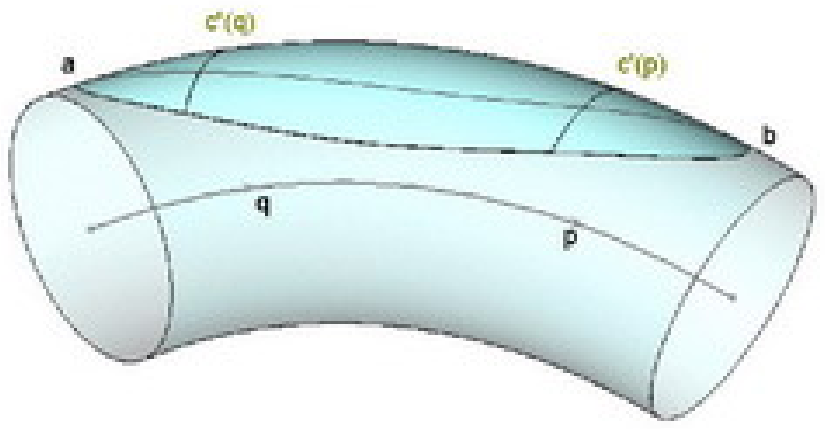}}

Again, $c'(p)$ is a continuous family of curves, here parametrized
by $p$. Since $\ii$ is parallel to the surface $\vv$, the surface
$S_{c'(p)}$ is contained in the normal plane of $p$. By the
curvature hypothesis, these plane sections intersect, say in
$\widetilde{\mathcal{I}}$, and by lemma~\ref{strates},
$\widetilde{\mathcal{I}}$ is either a point, either a straight line.
The later is impossible since, by lemma~\ref{etoile} the line
segment $\widetilde{\mathcal{I}}$ would be parallel to every $c'(p)$
whose curvatures are not zero. Therefore, $\widetilde{\mathcal{I}}$
has to be a point $q$, and all the normals of $\vv$ intersect in
$q$. Necessarily, the curve $\ii$ need to restrict to the point $q$.
In that case $\vv$ is parallel to $q$, hence $\vv$ is contained in a
sphere.

 \centerline{
 \epsfxsize=1,8in \epsfysize=1.3in \epsfbox{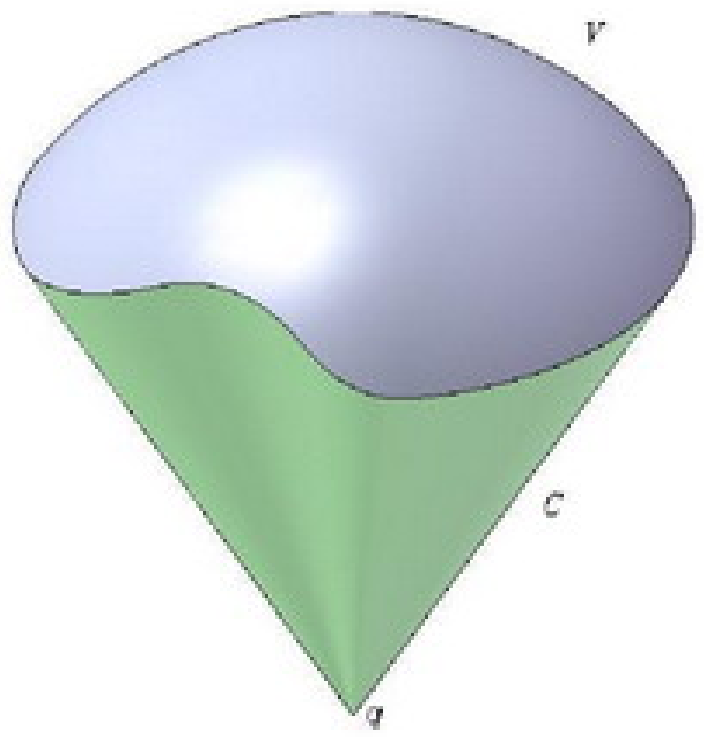}}

\noindent The Gaussian curvature on $\vv$ has to be constant, say
$\frac{1}{R^2}$, and all the points of $\vv$ are at distance $R$
from $q$. We are left to show that if we extend $\vv$ to the entire
leaf $\mm$, the distance between
$\mm$ and $q$ will be preserved, i.e. $\mm$ is a sphere.\\

Without loss of generality we consider $\vv$ to be the maximal
spherical cap with pole $y$. We denote $\mathcal{C}$ the cone with
apex $q$ generated by the normals of $\vv$, $\partial \vv$ the
boundary of $\vv$, and $\vv_{x}$ the curve obtained by extending
$\vv$ through $x\in
\partial \vv$ perpendicularly to $\partial \vv$ along the leaf. We consider as the $Z$-axis the line containing
 $q$ and $y$  and we define $\pp_{w}$, the \emph{\alignment plane},
containing the $Z$-axis and the point $w$.

 \centerline{
 \epsfxsize=2,4in \epsfysize=2,2in \epsfbox{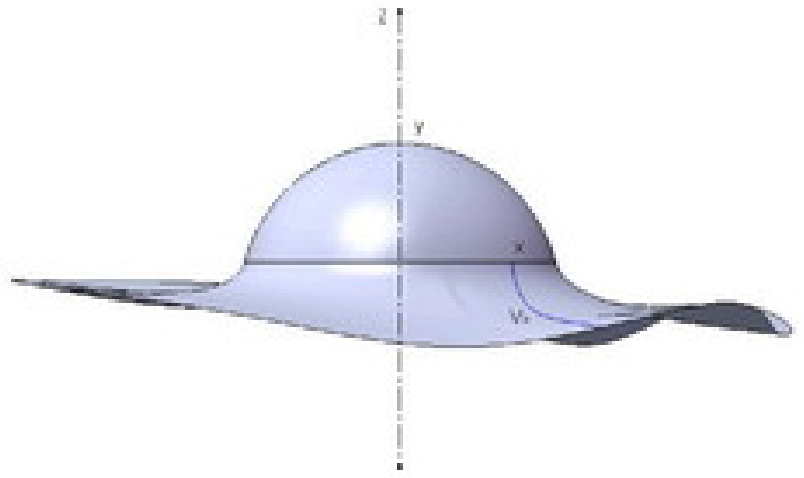}}
Remark that if the Gaussian curvature of $\mm$ changes along the
curve $\vv_{x}$, by
  smoothness of the leaves, for $y\in \partial \vv$ close to $x$,
  the curvature will also change along the curves $\vv_{y}$. The key
  point here is that if there is changes in the curvature,
   the surfaces $S_{\vv{x}}$ have to eventually leave
  their \alignment planes. Indeed we need to avoid two dimensional
  intersection with $\cc$, otherwise, together with the
  intersections
  of the
  surfaces $S_{\vv{y}}$ with $\cc$, we would get a three dimensional degenerate
  set. We are left with two possibilities. Either $\cc \cap
  S_{\vv_{x}}$ is always $q$, which is impossible since the curvature changes,
  either $S_{\vv_{x}}$ eventually leave the cone, which forces the normals of $S_{\vv{x}}$ to leave
  their
  \alignment plane. Note also that,
   if a normal stay in the alignment plane, it as to intersect $\cc$ at $q$.
   The main objective now is to show that if we extend the spherical cap, the normal lines
   stay in their alignment planes, intersecting the cone at $q$ and preserving the curvature of the spherical cap.\\

Let $\gamma_{\epsilon}$ be the closed curve in $\mm$ which is at
distance $\epsilon$ outside $\partial \vv$. By completeness, such
curve always exists for $\epsilon$ sufficiently small, say $\epsilon
< \varepsilon$. We want to show first that along such a path, all
the normal lines swing in the same direction with respect to their
\alignment planes. Let $x\in \gamma_\epsilon $, and assume the curve
is traversed in the clockwise direction.  Note that, if the dot
product between the normal line and the tangent vector of the curve
$\gamma_{\epsilon}$
 is positive, then there is an
increase of the $Z$-value of $\gamma_\epsilon$  around $x$. Suppose
we can take two curves, $\vv_x$ and $\vv_y$, for which the normal
lines swing in different directions. Say, without loss of
generality, that along each curve $\gamma_{\epsilon}$, $\epsilon \in
[ 0, \epsilon_1]$, either the direction changes only once, either
the normal is in the \alignment plane at $\vv_x$ and after
 roll in at most one direction. By smoothness of the leaf,
for each $\gamma_\epsilon$, there would be at least a point
$\alpha_\epsilon$ for which its normal lies in the \alignment plane
$\pp_{\alpha_\epsilon}$, hence intersects $\cc$ at $q$.  With an
appropriate choice of $\alpha_\epsilon$, for every $\epsilon$ in $[
0, \epsilon_1]$, we could generate the curve $\alpha(\epsilon)$
parallel to $q$. Since $\alpha(0) \in
\partial \vv$ the curve
 would be at distance $R$ of $q$.\\

Necessarily, there should be an other change of direction, say
between $\vv_y$ and $\vv_z$, along the curves $\gamma_\epsilon$ for
$\epsilon \in [0,\epsilon_2]$. We would then get an other curve
$\beta(\epsilon)$ at distance $R$ of $q$. Hence for all $\epsilon
\in [ 0, \tilde{\epsilon} ]$, where $\tilde{\epsilon}=\min \{
\epsilon_1, \epsilon_2 \}$, $\alpha_\epsilon$ and $\beta_\epsilon$
would be at distance $\epsilon$ from $\partial \vv$ and at distance
$R$ from $q$. Thus, for a fixed $\epsilon$, they would necessarily
have the same $Z$-value. But along $\gamma_\epsilon$, everywhere in
between $\alpha(\epsilon)$ and $\beta(\epsilon)$, the normal lines
are on the same side of their \alignment plane $\pp_z$, implying a
strict increase (or decrease) of the $Z$-value between the two
points. This is impossible, hence the
normal can only swing in one direction.\\

Therefore, the $Z$-value is monotonic as we follow the close curve
$\gamma_\epsilon$ in a given direction. This is impossible except if
the $Z$-value is constant, that is, if the normals stay in their
\alignment planes. Thus, for $\epsilon< \varepsilon$ the normal
lines of $\mm$ along $\gamma_\epsilon$ must intersect the cone at
$q$. Note that the curves $\vv_x$ stay parallel to $q$ when they
intersect $\gamma_\epsilon$. So we can increase $\vv$ to
$$\widetilde{\vv}:=\vv\bigcup_{\epsilon <
\varepsilon}\gamma_\epsilon$$ a bigger spherical cap. This
contradicts the maximality of $\vv$.
Hence $\vv$ has to be a sphere and is indeed the entire $\mm$. \hspace{\stretch{1}} $\Box$\\

\begin{Pro}Let $\mm$ be a leaf of the foliation $\Lambda$,
if there exists an open $\uu \subset \mm$ for which one of the
principal curvatures is identically zero and the other is
non-vanishing, then $\mm$ is an infinite cylinder.
\end{Pro}
\noindent{\sc {\bf Proof: }} Let $\lk(x)$ and ${\lk}'(x)$ be the
principal curves passing through $x\in \uu$ related respectively to
the vanishing  and the non-vanishing principal curvatures, say
$0~\equiv~\lambda_1~<~\lambda_2$. Note that $\lk(x)$ is a line.
Since $\lambda_2$ is never vanishing on $\uu$, the normal surfaces
$S_{\lk(x)}$ intersect, and from Lemmas \ref{strates} and
\ref{etoile}, all the lines $\lk(x)$ are parallel to $\ii$, which
has to be a line also.  By following the leaf along ${\lk}'(x)$, the
distance between $\ii$ and the points in $\uu$, say $R$, has to be
preserved. Therefore,
 $\uu$ has to be contained in a radius $R$ cylinder.

\centerline{
 \epsfxsize=2,5in \epsfysize=1,8in \epsfbox{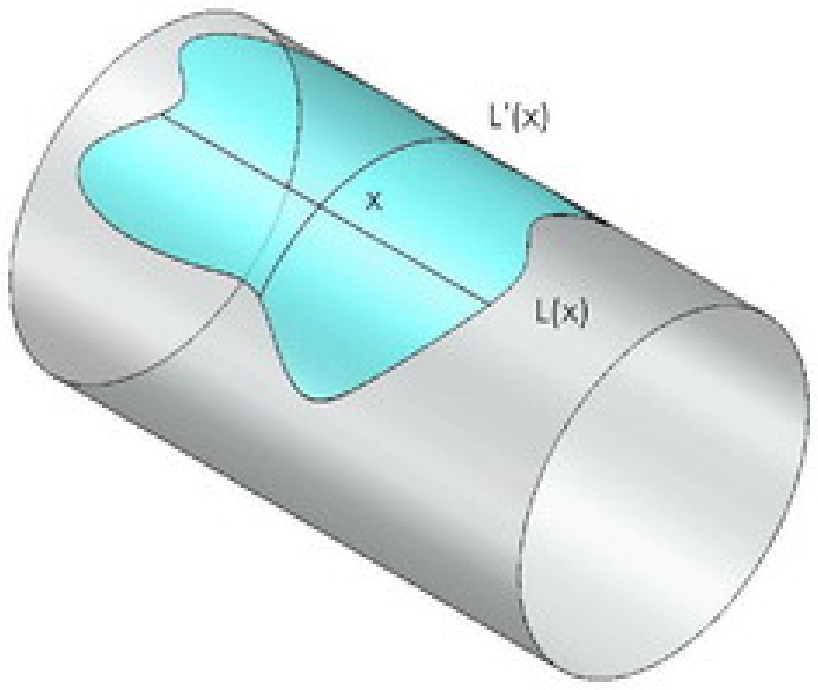}}

We can extend ${\lk}'(x)$ outside $\uu$, and for $\epsilon$
sufficiently small, we define the curve $\gamma_\epsilon$ to be the
union of the points along these extensions that are at distance
$\epsilon$ from the boundary of $\uu$. The normal surfaces
$S_{\gamma_\epsilon}$ have to intersect otherwise it would create a
three dimensional degenerate set with the normal lines of $\uu$. By
Lemma~\ref{strates}, these normal surfaces need to intersect in a
curve and since the leaf is smooth, this curve has to be
$\mathcal{I}$. By Lemma~\ref{etoile}, these curves are parallel to
the line $\ii$ so, by extending ${\lk}'(x)$ along $\mm$, we get a
cylinder $\cc$. \\

If we extend the principal curve $\lk(x)$ outside the cylinder
$\cc$, the curve obtained, say $\lk(x)^*$, needs to be a straight
line. Otherwise, as for the previous case, it would create a three
dimensional degenerate set while intersecting the normal lines of
the cylinder. Consequently, $\mm$ has to be an infinite cylinder.\hspace{\stretch{1}} $\Box$\\

\begin{Pro}If $\mm$ is a leaf of the foliation $\Lambda$,
then there is no open $\uu \subset \mm$ for which $\lambda_1\cdot
\lambda_2<0$.
\end{Pro}
\noindent{\sc {\bf Proof: }} Assume the opposite, and pick $y \in
\uu$. We consider $c(0)$, the intersection of $\uu$ with the
principal curve through $y$ associated to $\lambda_1<0$. We denote
$c(t)$, the curves of $\uu$ parallel to $c(0)$. By the curvature
hypothesis,  the normal surfaces $S_{c(t)}$, intersect. Hence by
Lemma~\ref{strates} and Lemma~\ref{etoile}, they
 intersect at $\mathcal{I}$ a line parallel to the curves $c(t)$.

\centerline{ \epsfxsize=2,4in \epsfysize=2in
\epsfbox{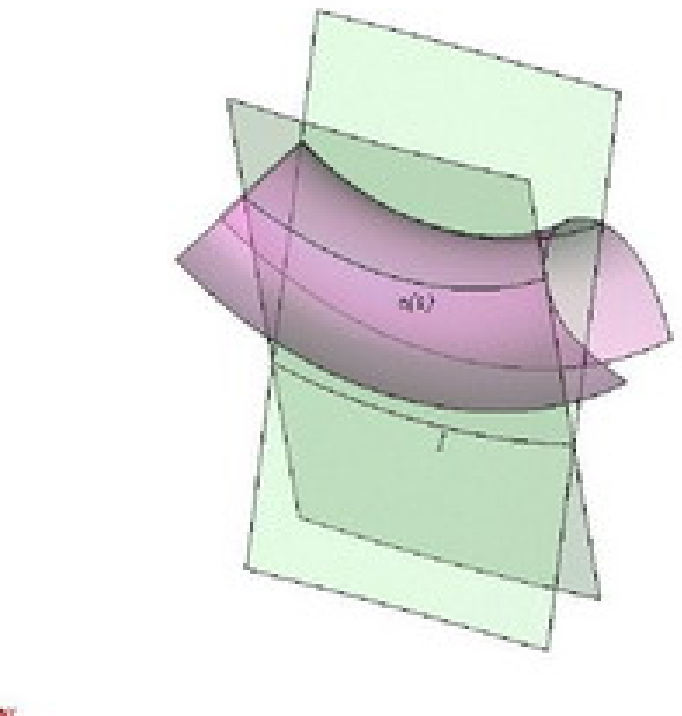}}

 After crossing the intersection curve, the normal lines of $\uu$ will cross
 leaves for which the two principal curvatures have the same sign,
  hence by Proposition~\ref{sphere}, the leaves on the other side of $\mathcal{I}$ will be
spheres. Such cohabitation of
 hyperbolic surface and spheres is impossible for the foliation $\Lambda$. \hspace{\stretch{1}} $\Box$\\

\begin{Co}\label{foliation}The non-degenerate leaves of the foliation
$\Lambda$ are either planes, cylinders or spheres.
\end{Co}

There seems to be a deep connection between the leaves arising from
a \imprimitive action and \emph{isoparametric manifolds}. Recall
that a hypersurface $\mathbf{M}^n$ of a Riemannian manifold
$\mathbf{V}^{n+1}$ is an isoparametric manifold if $\mathbf{M}^n$ is
locally a regular level set of a function $\lambda$ with the
property that both $\| \grad (\lambda)\|$ and $\Delta (\lambda)$ are
constant on the level sets of $\lambda$. One easily check that the
three possible leaves obtained in this section are indeed
isoprametric manifolds. The interesting point is that the only
complete isoprametric hypersurfaces of $\mathbb{R}^{3}$ are planes,
spheres and round cylinders and this classification holds for any
hypersurfaces of $\mathbb{R}^{n+1}$, see \cite{segre} for details.
Hence, the theory of isoparametric manifolds should provide a good
setting to approach Turbiner's conjecture in higher dimensions.

\section{$3$D-Turbiner's conjecture}

We have now all the tools needed to prove the $3$D modified
Turbiner's conjecture. This main theorem is a partial affirmation of
Turbiner's conjecture in three dimensions. First note that the
original conjecture involved complete separation while here we
succeed to prove that the equations separate partially. By a partial
separation,
 we mean that the equations separate into two equations, one involving only one
variable, the other involving the remaining variables. Also, three
assumptions need to be added to the original conjecture:  the
underlying action has to be \imprimitive \!\!, the manifold on which
the operator is defined has to be compact or can be compactified and
the contravariant metric need to be diagonal and \avoir\!. As for
the two dimensional case, the recipe is to pull back the invariant
foliation to the Euclidean environment where the leaves can only be
prescribed surfaces. Then working out the formulas for the operators
in the appropriate coordinates, one succeed to isolate one of the
variables. The major differences with the two dimensional case are:
the necessity of the \avoir requirement and the partial separation
obtained. As mentioned previously, at least one of the extra
requirements of the modified version of the conjecture is necessary;
a counter-example will be given in the next section.

\begin{Theo}\label{Turbiner}[3$D$ Modified Turbiner's conjecture]
Let $\hh_0$ be a second-order Lie algebraic operator generated by
the operators $T_a$ as per $(\ref{lie-alge})$, $g^{(ij)}$ be the
induced contravariant metric and $\mathbf{R}$ be a connected
component of $\mathbf{M_0}$ for which $g^{(ij)}$ is positive
definite.
 Suppose the following statements are true:
\begin{enumerate}
\item $\hh_0$ is gauge equivalent to a Schrödinger operator;
\item $(\mathbf{R}, g^{(ij)})$ is flat;
\item the operators $\{ T_a : a \in \g \}$ act \imprimitively;
\item $\mathbf{R}$ is either compact, or can be compactified in such a way that the $\mathbf{G}$-action on
$\mathbf{R}$ extends to a real-analytic action on the
compactification;
\item the metric $g^{(ij)}$ is diagonalizable and \avoir or $\mathbf{M}$ is a transverse, type changing manifold.
\end{enumerate}
Then, both the eigenvalue equation $\hh_0 \psi=E\psi$, and the
corresponding Schrödinger
equation separate partially in either a Cartesian, cylindrical or spherical coordinate system.\\
\end{Theo}
\noindent{\sc {\bf Proof: }} We denote $\Lambda$, the
$T_a$-invariant foliation.  The leafs are the level sets of a
function, say $\lambda$ and, from Proposition (\ref{imp}), this
foliation is also $\mathbf{G}$-invariant. The almost-Riemannian
manifold $(\mathbf{R}, g^{(ij)})$ fulfill the hypothesis of the
Tiling theorem,  thus there exist a real analytic map
$\Phi:\mathbb{R}^3\rightarrow \mathbf{R}$ for which $g^{(ij)}$ is
the push forward of the Euclidean metric. It is then possible to
pull back the rank $2$ foliation $\Lambda$ to get $\Phi^*(\Lambda)$
which is of rank two almost everywhere. From Corollary
\ref{perpendicular}, $\Phi^*(\Lambda)$ is locally orthogonal to a
foliation by geodesics that are, in this context, straight lines.
The rank two leaves are complete, hence, we can apply Corollary
\ref{foliation}, to conclude that there exists Cartesian coordinates
$(x,y,z)$ such that the leaves are given by the level sets of
$\lambda$, where
$\lambda$ is either $x$, $x^2+y^2$ or $x^2+y^2+z^2$.\\

We will now move the setting to $\mathbb{R}^3$. There is still the
local action of the group $\mathbf{G}$, but this action is
non-degenerate only whenever the Jacobian of $\Phi$ is not
degenerate. Separation is a local phenomenon, so for the present
purpose we can safely ignore the point of degeneracy\\

The operator $\hh_0$ is gauge equivalent to a Schrödinger operator
$\hh$, hence it must satisfy the closure condition. That means that
there exist a function $\sigma$ such that
\[ \hh_0=\Delta+\grad(\sigma)+V_0. \] From the \imprimitivity of the
action, $\hh_0(\lambda)=f(\lambda)$ and one easily verify that the
Laplacian of $\lambda$ is a function of $\lambda$ for the three
possible coordinate systems. Thus $\Lambda$ is also invariant with
respect to $\grad(\sigma)+V_0$. But, remark that
\[ [\grad(\sigma)+V_0](\lambda^2)-\lambda[\grad(\sigma)+V_0](\lambda)=\lambda
\grad(\sigma)(\lambda),\] which forces both $\grad(\sigma)$ and
$V_0$ to be function of $\lambda$. Depending of the metric, one
easily check that this force the gauge factor to separate the
following way:
\begin{eqnarray*}
\sigma(x,y,z)&=&\rho(x)+\eta(y,z),\\
\sigma(r,\theta,z)&=&\rho(r)+\eta(\theta,z), \\
\sigma(r,\theta,\phi)&=&\rho(r)+\eta(\theta,\phi).\end{eqnarray*}

Therefore the equation $\hh_0 \psi=E\psi$ separate partially and we
are left to show that the Schrödinger equation  also separate.\\

Recall that $V$, the potential of the Schrödinger operator is given by:\\
\[V=V_0+\grad(\sigma)^2+\Delta(\sigma), \textrm{ where }
V_0=V_0(\lambda).\]  After easy computations, the potentials are
given respectively by
\begin{eqnarray*}
V&=&F(x)+G(y,z),\\
V&=&F(r)+\frac{1}{r^2}G(\theta,z)+H(\theta,z),\\
V&=&F(r)+\frac{1}{r^2}G(\theta,\phi),
\end{eqnarray*}
where $F$ depends on $\rho$ and $V_0$, while $G$ and $H$ depend on
$\eta$.\\

This is sufficient to conclude that the Schrödinger equation \[
(\Delta+V) \Psi=E \Psi
\] separates partially either in Cartesian, cylindrical and
spherical coordinates. Indeed, we can perform respectively the
following separations:

\begin{eqnarray*}
[\partial_{xx}+F(x)-E]\Psi_1(x)&=&\alpha \Psi_1(x)\\
 \ [ \partial_{yy}+\partial_{zz}+G(y,z)]\Psi_2(y,z) &=&-\alpha
\Psi_2(y,z),\\
&&\\
 \ [
\partial_{rr}+\frac{1}{r}\partial_r+F(r)-E]\Psi_1(r)&=&(\frac{1}{r^2}\alpha+\beta)
\Psi_1(r)\\
   \ \  [  \partial_{\theta \theta}+\partial_{zz}+G(\theta,z)+H(\theta,z)]\Psi_2(\theta,z) &=&-(\alpha+\beta)
   \Psi_2(\theta,z),\\
&&\\
\  [
\partial_{rr}+\frac{2}{r}\partial_r+F(r)-E]\Psi_1(r)&=&\frac{1}{r^2}\alpha
\Psi_1(r)\\
   \ \  [ \frac{1}{\sin^2 \phi} \partial_{\theta \theta}+\partial_{\phi \phi}+\cot \phi \partial_\phi+
G(\theta, \phi) ]&=&-\alpha \Psi_2(\theta,\phi),
\end{eqnarray*}
where $\alpha$ and $\beta$ are separation
constants.~\hspace{\stretch{1}} $\Box$\\

\section{Counter-example} To conclude, we exhibit an example to
show that the extra hypotheses can not be omitted. Indeed, we
construct a Lie algebraic Schrödinger operator using generating
operators that act in a primitive way and we check that the
potential can not be separated, even partially. This counter-example
is the natural generalization of the one given in
\cite{Imprimitively} for the two dimensional case. It also motivates
the notion of almost-Riemannian manifold by realizing the quotient
of Euclidean space by an infinite reflection group. The general idea
for this type of construction is to find a set of basic invariants
and use them as coordinates.\\

This construction is a bit different from the usual one. Instead of
choosing first the coefficients to generate a Lie algebraic operator
and then verify the closure condition afterward, we proceed in an
different order. We first create an almost-Riemannian manifold
intimately related to the Lie algebra, then create an operator
satisfying the closure condition and finally check if there is a
choice of coefficients that generate that operator. We consider in
this example the Lie algebra $\mathfrak{sl}_4$, $\mathfrak{h}$ its
diagonal Cartan subalgebra equipped with the usual Killing inner
product and $W$, the affine Weil group associated to the root
system. We denote $L_1$, $L_2$, $L_3$ and $L_4$ the weights
associated to the diagonal entries of a trace-free diagonal matrix,
where $L_4=-L_1-L_2-L_3$. Taking $L_1$, $L_2$ and $L_3$ as
non-orthogonal coordinates, the contravariant form of the metric
tensor is given, in an appropriate basis, by

\[ \left( \begin{array}{ccc}
  2& -2/3 & -2/3 \\
 -2/3 & 2 & -2/3 \\
 -2/3 & -2/3 & 2\end{array} \right). \]

We define $z_k=e^{2\pi i L_k}$, the generators of the corresponding
 torus of diagonal unimodular matrices. The algebra of $W$-invariant
 elements of the complexified coordinate ring is generated by $\chi_1$, $\chi_2$ and $\chi_3$, the
 characters of the three fundamental representation of
 $\mathfrak{sl}_4\mathbb{C}$, see \cite{fulton} for more details.
 These three invariants are given by
\begin{eqnarray*}\chi_1&=&z_1+z_2+z_3+z_4,\\
  \chi_2&=&z_1 z_2+z_1 z_3+z_1 z_4+z_2 z_3+z_2 z_4+z_3 z_4,\\
 \chi_3&=&z_1 z_2 z_3+ z_1 z_2 z_4+ z_1 z_3 z_4+z_2 z_3 z_4,
  \end{eqnarray*}
and one easily compute the contravariant metric associated to this
algebra:
\[g^{(ij)}=-8\pi^2
\left(\begin{array}{ccc}
  \chi_1^2-8/3\chi_2
  &
 \frac{2}{3}(\chi_1
\chi_2-6\chi_3) &
\frac{1}{3}(\chi_1 \chi_3-16)\\
  \frac{2}{3}(\chi_1
\chi_2-6\chi_3)&\frac{4}{3}(\chi_2^2-2\chi_1 \chi_3-4) &
\frac{2}{3}(\chi_2 \chi_3-6\chi_1) \\
\frac{1}{3}(\chi_1 \chi_3-16) & \frac{2}{3}(\chi_2 \chi_3-6\chi_1) & \chi_3^2-8/3\chi_2\\
 \end{array}\right).\]

On the real torus, $\chi_1$ and $\chi_3$ are complex conjugates,
while $\chi_2$ is real-valued. Thus, fundamental invariants, denoted
$(x,y,z)$, are given by the real and imaginary parts of $\chi_1$ and
by $\chi_2$.
In the real coordinates, the corresponding contravariant metric
$g^{(ij)}$ , modulo a factor $\frac{-8 \pi^2}{3}$, reads as follow:
 \begin{eqnarray}\label{metrictensor}
\left(
\begin{array}{ccc}
  2x^2-z^2-4y-8& 2(xy-6y) & 3xz \\
2(xy-6y) & 4(y^2-2x^2-2z^2-4) & 2(yz+6z) \\
3xz & 2(yz+6z) & 2z^2-x^2+4y-8\end{array} \right).
\end{eqnarray}\\

For convenience, we will omit this $-8 \pi^2/3$ factor and one can
verifies that the Riemannian curvature
 tensor is identically zero where the metric is positive definite. The locus of degeneracy of the metric is given by

 \begin{eqnarray*}\sigma&=&-16(x^2+z^2)^3+(x^2+z^2+58/39)(320y^2+768)\\
 &&+(x^2-z^2)(32y^3-1152y)
+(x^4+z^4-352/39)(-4y^2+240)\\&&-144(x^4-z^4)y-8x^2y^2z^2-1248x^2z^2-64y^4=0.
\end{eqnarray*}

The objective now is to construct a Lie algebraic Schrödinger
operator on a space for which the contravariant metric is given by
(\ref{metrictensor}). The entries of the matrix are degree two
polynomials, hence the metric tensor can be generated by
$\mathfrak{a}_3$, the Lie algebra of infinitesimal affine
transformations of $\mathbb{R}^3$. A set of generators of
$\mathfrak{a}_3$ is given by:

\[ \begin{array}{cccccc} T_1=\partial_x, & T_2=\partial_y, & T_3=\partial_z, & T_4=x\partial_x, &
T_5=x\partial_y, & T_6=x\partial_z, \\ T_7=y\partial_x, &
T_8=y\partial_y, & T_9=y\partial_z, & T_{10}=z\partial_x, &
T_{11}=z\partial_y, & T_{12}=z\partial_z, \end{array} \] \noindent
and one easily sees that there is no function $\lambda$ for which
$\lambda$ and $T_\alpha(\lambda)$ are functionally dependant for all
the generators $T_\alpha$. Thus, these operators do not admit an
invariant foliation and the realization of $\mathfrak{a}_3$ is
therefore primitive. In term of these operators, the Laplacian, in
the $(x,y,z)$ coordinates, is given by
\begin{eqnarray*}
\Delta&=&-8T_1^2-16T_2^2-8T_3^2+2T_4^2-8T_5^2-T_6^2+4T_8^2-T_{10}^2-8T_{11}^2\\
  &
  &-2\{ T_1,T_7\}+2\{T_3,T_9\}-12\{T_2,T_4\}+12\{T_2,T_{12}\}+3\{T_4,T_{12}\}\\
  & & +2\{T_4,T_8\}+2\{T_8,T_{12}\}-2T_4-4T_8-3T_{12}.
\end{eqnarray*}
For $\sigma$ the determinant of the contravariant metric
$(\ref{metrictensor})$, one
 easily verifies that
\begin{equation}
\nabla \log \sigma=12(T_4+T_{12})+16T_8.
\end{equation}\\

Therefore, the operator
\[ \mathcal{H}_0=-\Delta+\nabla \log \sigma \]
 is Lie algebraic and gauge equivalent to $\hh$, a Schrödinger
operator, via the gauge transformation:

\[\hh=e^{-\log(\sigma)/2}\circ \hh_0 \circ
e^{\log(\sigma)/2}=-\Delta+U.\]  The potential $U$ can be computed,
\begin{eqnarray*}U&=&80-64[x^6+3x^4z^2+3x^2z^4+z^6-18x^4y+2x^2y^3-2y^3z^2+18yz^4\\&&+60x^4+60x^2y^2-312x^2z^2-8y^4+60y^2z^2+60z^4+\\&& -360x^2y+360yz^2+336x^2+192y^2+336z^2-640]\sigma^{-1},
\end{eqnarray*}
and can also be described in term of the affine coordinates $(L_1,
L_2, L_3)$ by
\begin{equation}\label{pot}U=80+\sum_{1\leq j<k\leq 4}\frac{1}{\sin^2(\pi i(
L_j-L_k))}.
\end{equation}

To conclude our counter-example, we need to show that the
Schrödinger equation $\hh$ cannot be solved, even partially, by
separation of variables. The potential here is symmetrical in the 3
variables, hence the separation in respect to one variable would
imply a separation in the other ones, thus a complete separation of
variables. The Schrödinger equation can be solved by separation of
variables in only eleven coordinate systems, nine of which (with the
exception of paraboloidal coordinates) are particular cases of the
 ellipsoidal coordinates. According to \cite{sep}, these
coordinates are: rectangular (Cartesian), circular cylinder,
elliptic cylinder, parabolic cylinder, spherical, conical,
parabolic, prolate spheroidal, oblate
spheroidal, paraboloidal, ellipsoidal coordinates.\\

An appropriate change of coordinates gives the orthonormal system
$(y_1,y_2,y_3)$ and one gets the following similar potential:
\begin{equation}U=80+\sum_{1\leq j<k\leq 3}\frac{1}{\sin^2(2\sqrt{2/3}\pi i(y_j\pm y_k)}
\end{equation}
Since the nine first coordinate systems are particular cases of the
last one, we only have to show that there is no separation possible
in the two last coordinate systems:
ellipsoidal and paraboloidal.\\

The ellipsoidal system of coordinates $(\xi_1, \xi_2, \xi_3)$ is
related to the Cartesian one by
\begin{eqnarray*}y_1&=&\sqrt{\frac{(\xi_1^2-a^2)(\xi_2^2-a^2)(\xi_3^2-a^2)}{a^2(a^2-b^2)}},\\
y_2&=&\sqrt{\frac{(\xi_1^2-b^2)(\xi_2^2-b^2)(\xi_3^2-b^2)}{b^2(b^2-a^2)}},\\
y_3&=&\frac{\xi_1 \xi_2 \xi_3}{ab}, \ \textrm{ where } \ \xi_1^2
\geq a^2  \geq \xi_2^2 \geq b^2 \geq \xi_3^2 \geq 0, \ \ \ \ \ \ \ \
\ \ \ \ \ \ \ \ \ \ \ \ \ \ \ \ \ \ \ \
\end{eqnarray*}
while, for the paraboloidal, we have:
\begin{eqnarray*}y_1&=&\sqrt{\frac{(\xi_1^2-a^2)(\xi_2^2-a^2)(\xi_3^2-a^2)}{(a^2-b^2)}},\\
y_2&=&\sqrt{\frac{(\xi_1^2-b^2)(\xi_2^2-b^2)(\xi_3^2-b^2)}{(b^2-a^2)}},\\
y_3&=&\sqrt{\frac{1}{2}(\xi_1^2+ \xi_2^2+ \xi_3^2-a^2-b^2)},\
\textrm{ where }\  \xi_1^2 \geq a^2  \geq \xi_2^2 \geq b^2 \geq
\xi_3^2 \geq 0.
\end{eqnarray*}
For these two systems, a given potential $U$ separate if and only if
it is of the form
\[ U=\frac{(\xi_2^2-\xi_3^2)U_1(\xi_1)+(\xi_1^2-\xi_3^2)U_2(\xi_2)+(\xi_1^2-\xi_2^2)U_3(\xi_3)}
{(\xi_1^2-\xi_2^2)(\xi_2^2-\xi_3^2)(\xi_1^2-\xi_3^2)}. \]

After the suitable substitution, the potential $U$, in terms of the
new coordinates $(\xi_1, \xi_2, \xi_3)$, fails to be of that
requested form. Therefore,
no separation is possible.\\

Thus, this example emphasis on the necessity of the extra hypotheses
we needed to add to the original conjecture. Here, at least one of
these hypotheses, the imprimitivity of the action, fails to be
satisfied and the Schrödinger equation can not be solved by
separation of variables, even partially.\\


\subsection*{Acknowledgments:}
The research is supported by NSERC Grant $\#$RGPIN $105490-2004$ and
by a McGill Graduate Studies Fellowship. I would like to thank Rob
Milson for
   all the time he devoted to answer my many questions and for the feedback he gave me about this work.
   I would also like to thank Niky Kamran, for all
   these precious advice and all the encouragement he gave me.


\end{document}